\documentstyle[12pt]{article}
\newcounter{local}
\newcounter{locallocal}
\newcommand{\scl}{\stepcounter{local}}
\newcommand{\g}{\nabla \times  }
\newcommand{\vh}{ {\bf H}}

\newcommand{\z}{\varepsilon}
\newcommand{\h}{\hspace{1cm}}
\newcommand{\hh}{\hspace{2cm}}
\newcommand{\hhh}{\hspace{3cm}}

\newcommand{\by}{\begin{eqnarray}}
\newcommand{\ey}{\end{eqnarray}}
\newcommand{\bys}{\begin{eqnarray*}}
\newcommand{\eys}{\end{eqnarray*}}
\newcommand{\B}{ {\bf H}}
\newcommand{\E}{ {\bf E}}
\newcommand{\F}{ {\bf F}}
\newcommand{\J}{ {\bf J} }
\newcommand{\Ph}{ {\bf \Phi}}
\newcommand{\Ps}{ {\bf \Psi}}
\newcommand{\W}{{\bf W}}
\newcommand{\y}{{\bf c}}
\newcommand{\e}{{\bf e}}
\newcommand{\Y}{{\bf Y}}
\begin{document}
\begin{center}
\begin{bf}
On a Singular Limit Problem for Nonlinear Maxwell's Equations

\end{bf}
\ \\
\vspace{.5cm}
Hong-Ming Yin\\
Department of Mathematics,
University of Notre Dame\\
 Notre Dame, IN 46556.
\end{center}
\vspace{1cm}

\vspace{3cm}
{\bf Abstract}: 
In this paper we study 
the following nonlinear Maxwell's equations \\
$\z \E_{t}+\sigma(x,|\E|)\E= \g \vh +\F,\, \vh_{t}+\g \E=0$, where
$\sigma(x,s)$ is a monotone graph of $s$. It is shown that the system
has a unique weak solution. Moreover, 
the limit of the solution as $\z\rightarrow 0$ converges
to  the solution of quasi-stationary Maxwell's equations.

\ \\
{\bf AMS(MOS)} Subject Classifications: 35K20, 35Q20.

\ \\
{\bf Key Words and Phases}: Nonlinear Maxwell's Equations, Singular limit.

\newpage

\begin{center}
{\bf 1. Introduction}
\end{center} 
Let $\Omega$ be a bounded domain in $R^3$ and $Q_{T}=\Omega\times (0,T]$
for any fixed $T>0$. 
Let $\E$ and $\vh$ be the electric and magnetic fields, respectively,
in $\Omega$ (here and thereafter
a bold letter represents a vector in $R^3$).
Let $\sigma$ be the electric conductivity in the field, which
is assumed to be a function of $x$ and $|\E|$.
Consider the following Maxwell's equations (see
Landau-Lifschitz \cite{LL}):
\setcounter{section}{1}
\setcounter{local}{1}
\begin{eqnarray}
& & \z \E_{t}+\sigma(x, |\E|)\E=\g \vh + \F, \h (x,t)\in Q_{T},\\ 
& & \vh_{t}+\g\E= 0, \hhh (x,t)\in Q_{T},\scl\\
& & \nabla \times {\E}=0,
\hhh (x,t)\in \partial \Omega\times (0,T],\scl \\
& & \E(x,0)=\E_0(x), {\bf H}(x,0)={\bf H}_{0}(x), \h x\in \Omega, \scl
\end{eqnarray}
where $\z$ is the dielectric parameter and other physical parameters
are normalized.

 In some applications (\cite{DL,MM}), 
the electric conductivity, $\sigma$, strongly depends on the electric field
$|\E|$, hence the electric current density. Particularly, 
the electric conductivity may act like a switch-like function
in some electromagnetic fields. 
On the other hand, for many types of  micron devices  and
other industrial problems (such as microwave heating \cite{MM,GY,Y3}) etc.)
the experiment shows that
the displacement current, $\z \E_{t}$, is often negligible
since it is small in comparison of the eddy current, $\J=\sigma \E$. 
This motivates us to investigate the nonlinear 
problem (1.1)-(1.4) and the singular limit problem
as $\z\rightarrow 0$. 
It is shown that there exists a unique global solution
 to (1.1)-(1.4). Moreover, the limit of the solution converges to
 the solution of
the quasi-stationary system (i.e.,
the system (1.1)-(1.4) with $\z=0$ in (1.1)). This limit solution 
provides new existence result
for the quasi-stationary system. Indeed,
when $\z=0$, the system (1.1) becomes
\by 
\sigma(x, |\E|)\E =\g \B +\F,\scl
\ey
Thus, one can solve Eq. (1.5) for $|\E|$ in terms of $|\g \vh|$ and known
data,
\[ |\E|=g(x,|\g \B|),\]
where $g(x,s)$ is the inverse function of $\sigma(x,s)s$.

It follows from (1.2) that $\B$ satisfies
\by
\vh_t+\g [ \rho(x,|\g \vh|)\g\B ]=0,\scl
\ey
where
\[ \rho(x,|\g \vh||)=\frac{1}{\sigma(x, |\E|)}=
\frac{1}{\sigma(x, g(x,|\g\vh|))}\]
represents the electric resistivity in the field.

 The research on Maxwell's equations
  is of great interesting because of the important
   applications in plasma physics, 
    semiconductor-superconductor
     modeling and other industrial problems (\cite{GY,H,LL,Y3} etc.).
The study on the system (1.1)-(1.4) as well as the quasi-stationary form (1.6)
received considerable attention recently. In \cite{MP}, the authors established
the well-posedness  for
a quasi-stationary system, where a constitutive relation between
the magnetic field $\vh$ and the magnetic induction ${\bf B}$ is assumed to
be nonlinear.
In \cite{Y1}, the author studied
the regularity of weak solution to a linear system of (1.6) with
minimal requirement on coefficients. There is a special
interest  when
\[ \rho(|\g \vh|)=|\g \vh|^{p-2}, \h p>2.\]
On one hand, if $\B$ is restricted in one direction (scalar field) then
the evolution system (1.6) becomes the p-Laplacian which has been studied
extensively (see \cite{D} and the references therein). On the other hand,
in a recent work
\cite{Y2} (also see \cite{AEW,BP} for the scalar case),
it is shown that
the limit of the solution to (1.6)
as $p\rightarrow \infty$ is the unique solution to Bean's critical-state 
model
in the superconductivity theory (\cite{B}).
Thus, for large $p$ the system (1.6) provides a good approximation
to Bean's model.  More recently, 
 the author of \cite{JO} studied the similar problem to this paper
 in a domain with a bounded complement
  in $R^3$. The conditions on $\sigma$ in \cite{JO} is quite different from
   ours here. Like many nonlinear problems, the major difficulty is
   how to pass the weak limit of an  approximate solution for 
   a nonlinear function $\sigma(x,s)$. This is done by  employing a 
   monotonicity
   argument (\cite{E}).
   The monotonicity of $\sigma(x, s)$ in $s$ is essential 
   in the proof. 

 In \S 2, we use the finite element method to establish the well-posedness
 of the system (1.1)-(1.4) for fixed $\z>0$. 
 In \S 3, we show that the singular limit of
 the solution to (1.1)-(1.4) has a unique limit. Moreover, the limit solution
 solves the quasi-stationary Maxwell's equations. Some examples
 are also discussed in this section.
 
 
\begin{center}
{ \bf 2. Existence and Uniqueness for fixed $\z>0$ }
\end{center}

 Introduce some standard spaces (see \cite{DL,GR}).
 \bys
 & &
  H(div, \Omega)=\{ {\bf U} \in L^2(\Omega)^3, div {\bf U}\in L^2(\Omega)\};\\
  & & H(curl,\Omega)=\{ {\bf U}\in L^2(\Omega)^3, curl {\bf U}\in L^2(\Omega)^3\}\, ;\\
  & & H(div0, \Omega)=\{ {\bf U}\in H(div, \Omega): div {\bf U}=0
  \mbox{ in $\Omega$} \},\\
  & & H_0(curl, \Omega)=\{ {\bf U}\in H(curl, \Omega): {\bf N}\times {\bf U}=0\,
  \mbox{on $\partial \Omega$} \},
  \eys
  where ${\bf N}$ is the exterior unit normal on $\partial \Omega$.

  Note that the trace of a function in $H(curl, \Omega)$ is well defined
  (see \cite{DL} for example).

We shall assume the following conditions on $\sigma(x,s)$ and
data $\E_0(x), 
{\bf H}_{0}(x)$ and ${\bf F}(x,t)$.
\ \\
H(2.1): Let $\sigma(x, s)$ be measurable 
in $\Omega\times [0, \infty)$ and monotone increasing in $s$. Moreover,
\bys
& & \int_{0}^{s^{2} } \sigma(x,\sqrt{s}) dx\geq a_0 s^{p+2}-a_1, \mbox{
if $s$ is sufficiently large},\\
& & 0\leq \sigma(x,s)\leq b_0 (1+s^{p}), s\in [0,\infty),\mbox{for some $p\geq 
0$},
\eys
where the constants $a_0>0, a_1\geq 0$ and $b_0\geq 0$.

\ \\
H(2.2): Assume that
$ {\bf H}_{0}\in H(curl, \Omega)\bigcap H(div0,\Omega),
{\bf F}\in H^1(0,T; H^1(\Omega))$.

\ \\
{\bf Definition 2.1}: {\em A pair of vector fields ($\E(x,t), \B(x,t))$ 
is said to be a weak solution of the problem (1.1)-(1.4), if 
\[ \E\in L^{2}(0,T; H_0(curl, \Omega))\bigcap L^{p+2}(0,T;\Omega), 
{\bf H}\in L^2(0,T;H(div0,\Omega) )\bigcap H(0,T, L^2(\Omega))\]
which satisfy the following
integral identities:
\setcounter{section}{2}
\setcounter{local}{1}
\by
& & \int
\int_{Q_{T}}\left[-\z \E\cdot {\bf \Phi}_{t}+ \sigma(x,|\E|)\E\cdot
\Ph \right] dxdt\nonumber\\
& & = \int\int_{Q_{T}}{\bf H}\cdot (\g \Ph)dxdt
+\z \int_{\Omega} \E_0\cdot \Ph(x,0)dx,\\
& & \int\int_{Q_{T}}\left[- \B
\cdot \Ps_{t} + \E \cdot (\g \Ps)\right] dxdt=
\int_{\Omega} \left[ {\bf H}_{0}(x) \cdot \Ps(x,0)\right]dx \scl
\ey
for all test functions $\Ph\in H^1(0,T; H_0(curl, \Omega)), \Ps\in H^1(0,T; 
H(curl, \Omega))$ with \\$\Ph(x,T)=\Ps(x,T)=0$.}

First of all, we derive some energy estimates. 
A special attention is paid on how various constants depend on $
\z$ since we will study the singular limit problem in section 3.

\ \\
{\bf Lemma 2.1}: {\em Under the assumptions H(2.1)-H(2.2)
there exist constants $C_1, C_2$ and $C_3$  such that
\setcounter{local}{2}
\bys
& &\sup_{[0,T]}\int_{\Omega}\left[\z |\E|^2+ |{\bf H}|^2
+|\E|^{p+2}\right] dx
+ \int_{0}^{T}\int_{\Omega} \left[\z |\E_{t}|^2+|\vh_{t}|^2\right]dxdt
\nonumber\\
& & \leq C_1\int_{\Omega} \left[|\E_0|^2+ |{\bf H}_0|^2+ |\g \B_0|^2
\right]dx
+C_2\int\int_{Q_{T}}[|\F|^2+|\F_t|^2] dxdt+ C_3,
\eys
where $C_1, C_2$ and $C_3$ depend only on known data.}\\
{\bf Proof}: Note that for any vector fields ${\bf A}, {\bf B}\in H(curl, 
\Omega)$ with either ${\bf A}$ or ${\bf B}$ in $H_0(curl, \Omega)$,
the following identity holds:
\[
\int _{\Omega} {\bf A}\cdot ({\bf \g B}) dx=
\int _{\Omega} {\bf B}\cdot ({\bf \g A}) dx.\]
Taking inner product to the system (1.1) and (1.2) by ${\bf E}$ and $\B$, 
respectively, we add up the resulting equations 
to obtain 
\by
& &
\sup_{[0,T]}\int_{\Omega}\left[ \z |\E|^2+|{\bf H}|^2\right] dx
+
\int\int_{Q_{T}}\sigma(x, |\E|)
|\E|^{2}dxdt\nonumber \\
& & \leq C\int_{\Omega}\left[|\E_0|^2+{\bf H}_{0}|^2\right]dx 
+\int\int_{Q_{T}}[|\E\cdot\F|]dxdt\scl
\ey
where the constant $C$ depends only on known data, but not on $\z$.

We first assume that $\sigma(x,s)$ is differentiable with respect to $s$.
Then we formally differentiate Eq.(1.1) and Eq.(1.2) with respect to $t$ 
to obtain 
\bys
& &\z \E_{tt}+\sigma(x,|\E|)\E_t+\sigma_{s}(x,|\E|)(|\E|)_{t}\E=\g \B_t+\F_t,\\
& & \B_{tt}+\g\E_t=0.
\eys
It is clear that
\[ \int\int_{Q_{T}}(\g \E_t)\cdot \B_t dxdt=
\int\int_{Q_{T}} (\g\B_t)\cdot \E_tdxdt.\]
We take the inner product by $\E_t$ for the first equation and by
$\B_t$ for the second equation and add up the resulting equations
to obtain:
\bys
& & \sup_{[0,T]}\int_{\Omega}\left[ \z |\E_t|^2+|\B_t|^2\right]dxdt+
\int\int_{Q_{T}}\left[\sigma|\E_{t}|^2+\sigma_s (|\E|)_{t}\E\cdot \E_{t}\right]
dxdt\\
& & \leq C,
\eys
where $C$ depends only on known data.

Note that $\sigma_{s}\geq 0$, we see that
\bys
& & \int\int_{Q_{T}}\sigma_{s}(|\E|)_{t}\E\cdot \E_{t} dxdt\\
& & =\int\int_{Q_{T}}\sigma_{s}(|\E|)_{t} \frac{d}{dt}\frac{|\E|^{2}}{2}dxdt\\
& & =\int\int_{Q_{T}}\sigma_{s}|\E|[ |\E|_{t}]^2dxdt\geq 0.
\eys
It follows that
\bys
& & \sup_{[0,T]}\int_{\Omega}\left[ \z |\E_t|^2+|\B_t|^2\right]dxdt+
\int\int_{Q_{T}}[\sigma|\E_{t}|^2 dxdt
\leq C.
\eys
Since the above estimate does not depend on the differentiability of $\sigma$
with respect to $s$, therefore 
the above estimate holds as long as $\sigma$ is
monotone increasing with respect to $s$.

Now we take the inner product by $\E_t$ to (1.1) and by $\B_t$ to (1.2)
and then add up the resulting equations to obtain
\bys
& & \int\int_{Q_{T}}\left[ \z |\E_{t}|^2+|\vh_{t}|^2\right]dxdt+
\frac{1}{2}\int_{\Omega}\int_{0}^{|\E(x,T)|^{2}}\sigma(x, \sqrt{s})dsdx\\
& &\leq
\int\int_{Q_{T}}\E_{t}\cdot \F dxdt +C\int_{\Omega}\int_{0}^{|\E_{0}(x)|^{2}}\sigma(x, \sqrt{s})dsdx+\int_{\Omega} |\g \B_0|^2dx+C.
\eys
Now
\bys
& & \int\int_{Q_{T}}\E_{t}\cdot \F dxdt\\
& & =\int_{\Omega} \left[\E(x,T)\cdot \F(x,T)-\E_{0}(x)\cdot\F(x,0)\right]dx
-\int\int_{Q_{T}}\E\cdot \F_{t}dxdt\\
& & \leq \int_{\Omega}\left[\frac{a_{0}}{4}\E(x,T)|^{p+2}+\frac{16}{a_{0}}
|\F_{t}|^{\frac{p+2}{p+1}}\right] dxdt.
\eys
On the other hand,
by the assumption H(2.1) we may assume that the growth condition of 
$\sigma(x,s)$ on $s$
holds for all $s\geq M_0$, i.e.,
\[ \int_{0}^{s^{2} } \sigma(x,\sqrt{s}) dx\geq a_0 s^{p+2}-a_1,
\mbox{ if $s\geq M_0$},\]
where $M_0$ is a fixed constant.

It follows that
\by
& & \int_{\Omega}\int_{0}^{|\E(x,T)|^{2}}\sigma(x,s)dsdx\nonumber\\
& & \geq a_0 \int_{\Omega \bigcap \{x:|\E(x,T)|\geq M_{0}\}}
|\E(x,T)|^{p+2}dx -C.\scl
\ey
Combining (2.3)-(2.4) yields
\bys
& & \sup_{0\leq t\leq T}\int_{\Omega}[\z |\E|^2+|\vh|^2]dx+\sup_{0\leq t\leq T}
\int_{\Omega}|\E(x,t)|^{p+2}dx+
\int_{0}^{T}\int_{\Omega}[\z |\E_t|^2+|\vh_t|^2]dxdt\\
& & \leq C\int\int_{Q_{T}}
[|\F|^{2}+|\F_{t}|^{2} ]dxdt+\int_{\Omega}[|\E_0|^2+|\B_0|^2+|\g\B_0|^2]dx+C.
\eys

\hfill Q.E.D.

\ \\
{\bf Theorem 2.2}: {\em Under the assumptions H(2.1)-H(2.2) the
problem (1.1)-(1.4) has a unique weak solution. Moreover,
\[ curl \E\in L^2(Q_{T}), \E_t\in L^2(Q_{T})\]
and} 
\[  {\bf H}_{t}\in L^{2}(Q_{T}), \g\B\in L^{\frac{p+2}{p+1}}(Q_{T}).\]
{\bf Proof}: 
The proof is based on the finite element method (see \cite{LSU} for
parabolic equations).
The monotonicity of $\sigma(x, s)$ on $s$ plays
an important role. We shall first deal with the case where $\sigma(x, s)$
is continuous on $s$.  For convenience, we
rewrite the system (1.1)-(1.4) to the following form:
\by
& &\z \W_{tt}+\sigma(x,|\W_t|)\W_{t}=\g[ \B_0-\g \W] + \F, \h (x,t)\in Q_{T},
\scl \\
& &{\bf N}\times (\W_{t})=0, \hhh x\in \partial \Omega, 0\leq t\leq T, \scl \\
& & \W(x,0)=0, \W_t(x,0)=\E_0(x), \hh x\in \Omega,\scl
\ey
where $\W$ is defined as follows:
\[ \W(x,t)=\int_{0}^{t}\E(x,\tau)d\tau.\]
It is clear that if 
$\W$ is a solution of the system (2.5)-(2.7) then 
a pair of functions defined by
\[ \E(x,t)=\W_t(x,t), \B(x,t)=\B_0(x)-\g \W(x,t)\]
will be a weak solution of (1.1)-(1.4).
Let $\{ {\bf e}_{k}\}=\{ e_{k}^{(1)}, e_{k}^{(2)}, e_{k}^{(3)}\}
$ be a smooth basis of $H_0(curl, \Omega)$
and orthonormal in $L^2(\Omega)^3$, i.e.
\[ < {\bf e}_{i}, {\bf e}_{j}>=\delta_{ij},\]
where $\delta_{ij}=1$ if $i=j$ and $\delta_{ij}=0$ if $i\neq j$.

Now we expand the known data as follows:
\bys
& & \B_0(x)=\sum_{k=1}^{\infty}diag[{\bf a}_{k}\circ {\bf e}_{k}],\\
& & \E_{0}(x)=\sum_{k=1}^{\infty}diag[ {\bf b}_{k}\circ {\bf e}_{k}],\\
& & \F(x,t)=\sum_{k=1}^{\infty}diag[{\bf g}_{k}(t)\circ {\bf e}_{k}],
\eys
where ${\bf a}_{k}, {\bf b}_{k}$ and ${\bf g}_{k}$ are
$3\times 1$ matrices, the symbol $\circ$ is the matrix product and
$diag[ ]$ represents the diagonal vector of a matrix.

Let 
\[ \W_n(x,t)=\sum_{k=1}^{n}diag[ \y_{n}^{(k)}\circ {\bf e}_{k}],\]
where $\y_{n}^{(k)}(t)$ is a $3\times 1$ vector 
which is determined by the following ordinary differential
system:
\by
& & \z \frac{d^{2}}{dt^{2}}\y_{n}^{(k)}+ \sigma(x, |\W_{nt}|)
\frac{d}{dt} \y_{n}^{(k)}
=A_{k}(\W_{n}^{(k)}, {\bf e}_{k})+B_{k}(t),\scl \\
& & \y_{n}^{(k)})(0)=0,\scl \\
& & \frac{d}{dt}\y_{n}^{(k)}(0)={\bf b}_{k}, \scl
\ey
where
\bys
& & \W_{n}^{(k)}=diag[ \y_{n}^{(k)}\circ {\bf e}_{k}],\\
& & A_{k}(\W_{n}^{(k)}, {\bf e}_{k})=\int_{\Omega}diag\{
[\g \W_{n}^{(k)}]
\circ [\g {\bf e}_{k}]\}dx,\\
& & B_{k}(t) =\int_{\Omega}diag\{
(\g \B_0+\F)\circ {\bf e}_{k}\} dx, \h k=1,2,\cdots n.
\eys
Now we define the approximate solution $(\E_n, \B_{n})$ as follows:
\[ \E_{n}(x,t)=\W_{n t}(x,t), \B_n(x,t)=\B_{0n}-\g \W_{n}(x,t),\]
where 
\[ \B_{0n}(x)=\sum_{k=1}^{n}diag[{\bf a}_{k}\circ\e_{k}].\]
Equivalently, then $(\E_n, \B_n)$ satisfies the following system
in the weak sense:
\by
& & \z \E_{n t}+\sigma(x,|\E_{n}|)\E_{n}=\g \B_{n}, \h (x,t)\in Q_{T},\scl \\
& & \B_{n t}+\g \E_{n}=0, \hh (x,t)\in Q_{T}.\scl
\ey
Similar to Lemma 2.1, one can easily derive the following energy estimates
:
\bys
& &\sup_{[0,T]}\int_{\Omega}\left[\z |\E_n|^2+ |{\bf H}_n|^2
+|\E_n|^{p+2}\right]dx+
\int_{0}^{T}\int_{\Omega} \left[\z |\E_{nt}|^2+|\vh_{nt}|^2\right]dxdt
\nonumber\\
& & \leq C_1\int_{\Omega} \left[|\E_{0n}|^2+ |{\bf H}_{0n}|^2+ |\g \B_{0n}|^2
\right]dx
+C_2\int\int_{Q_{T}}[|\F_n|^2+|\F_{n t}|^2] dxdt+ C_3,
\eys
where $C_1,C_2$ and $C_3$ are independent of $n$ and $\z$.

By the weak compactness property, we can extract a subsequence (still denoted by
$(\E_{n}, \vh_{n}))$ such that
\bys
& & \E_{n}\rightarrow \E, \E_{n t}\rightarrow \E_t, \B_{nt}
\rightarrow \B_t, \, 
\mbox{weakly in $L^2(Q_{T})$,}\\
& & \vh_{n}\rightarrow \vh, \mbox{ weakly in $L^2(0,T; W^{1,\frac{p+2}{p+1}}
(\Omega))$},\\
& & \vh_n \rightarrow \vh, \mbox{ a.e. in $Q_{T}$}.
\eys
Moreover, 
\bys
\E_{n}\rightarrow \E \h \mbox{ weakly in $L^{p+2}(Q_{T})$}.
\eys

Next we claim that the sequence $\E_{n}$ converges to $\E$ strongly in
$L^2(Q_{T})^3$. To prove the claim we only need to show that $\{\E_{n}\}$ is
a Cauchy sequence in $L^2(Q_{T})^3$. Let 
\[ {\E}_{n}^*(x,t)= \E_{n}(x,t)-\E_{m}(x,t), \vh_{n}^*(x,t)=\vh_{n}(x,t)-\vh_{m}(x,t).\]
By energy estimates, we see
\bys
& & \sup_{0\leq t\leq T}\int_{\Omega}[ {|\E}_{n}^{*}|^2+|\vh_{n}^{*}|^2]dx
+\\
& & 
\int\int_{Q_{T}}\left\{[ \sigma(x, |\E_{n}|)\E_{n}-\sigma(x, |\E_{m}|)\E_{m}]
\cdot [\E_{n}-\E_{m}]\right\} dxdt\\
& & \leq C\int_{\Omega}[ |\E_{0n}-\E_{0m}|^2+|\vh_{0n}-\vh_{0m}|^2]dx+
C \int\int_{Q_{T}}[ |\F_n-\F_m|^2 ]dxdt,
\eys
where $C$ is a constant independent of $n$ and $m$.

Note that $\sigma(x,s)$ is monotonic increasing in $s$, then
\bys
& & [\sigma(x, |\E_{n}|)\E_{n}-\sigma(x, |\E_{m}|)\E_{m}]
\cdot [\E_{n}-\E_{m}]\\
& & \geq \frac{ \sigma(x, |\E_{n}|)-\sigma(x, |\E_{m}|)][|\E_{n}|^2-
|\E_m|^2]}{ 2}\geq 0.
\eys
It follows that
\bys
& & \sup_{0\leq t\leq T}\int_{\Omega}[ |{\E}_{n}^{*}|^2+|\vh_{n}^{*}|^2dx\\
& & \leq  C\int_{\Omega}[ |\E_{0n}-\E_{0m}|^2+|\vh_{0n}-\vh_{0m}|^2] dx+
C\int\int_{Q_{T}}|\F_n-\F_m|^2dxdt.
\eys
This implies that $ \E_{n}, \vh_{n}$ are Cauchy sequences since
both $\E_{0n}, \vh_{0n}$ and $\F_n$ are Cauchy sequences in $L^2(Q_{T})^3$.
Hence, 
\[ \E_n, \B_n\rightarrow \E, \B \, \mbox{strongly in $L^2(Q_{T})$}.\]
After taking a subsequence if necessary, we see that
\[ \E_{n}\rightarrow \E, \h a.e. \,in \, Q_{T}.\]

To show the existence of a weak solution to (1.1)-(1.4), we only need to show
\[ \sigma(x, |\E_n|)\E_n\rightarrow \sigma(x, |\E|) \E
\h \mbox{ in $L^1(Q_{T})$}.\]

As $\sigma(x,s)$ is continuous on $s$ and $\E_{n}$ converges to $\E$ almost everywhere
in $Q_{T}$, we know 
\[ \sigma(x, |\E_n|)\E_n\rightarrow \sigma(x, |\E|) \E
\h a.e. \mbox{ in $Q_{T}$.} \]
We now show that $\sigma(x, |\E_n|)\E_n$ is equip-integrable in $Q_{T}$.
We adopt a technique used for scalar elliptic and parabolic equations.
Let $A$ be any measurable subset of $Q_{T}$. For any large $m>0$,
\bys
& & \int\int_{A}\sigma(x, |\E_{n})|\E_{n}|dxdt\\
& & \leq \int\int_{A\bigcap \{|\E_{n}|\leq m\}}\sigma(x, |\E_{n}|)|\E_n|dxdt+
\int\int_{A\bigcap \{|\E_{n}|\geq m\}} \sigma(x, |\E_{n}|)|\E_n|dxdt\\
& & \equiv I_1+ I_2.
\eys
The assumption on $\sigma(x,s)$ yields 
\bys
& & I_1\leq C\int\int_{A\bigcap \{|\E_{n}|\leq m\}}[1+|\E_{n}|^{p}
|\E_n|]dxdt,
\eys
which can be arbitrarily small if $|A|$ is small since 
$\E_{n}\in L^{p+2}(Q_{T})$.

On the other hand,
\bys
& & I_2 \leq \frac{1}{m} \int\int_{A\bigcap \{|\E_{n}|\geq m\}}
\sigma(x, |\E_{n}|)|\E_{n}|^2dxdt\leq \frac{C}{m},
\eys
which is also small if $ m$ is sufficiently large.

This concludes that $\sigma(x, |\E_n|) \E_{n}$ is equip-integrable in $Q_{T}$.

It follows by Vitali's theorem that
\[ \sigma(x, |\E_n|) \E_{n}\rightarrow \sigma(x, |\E|)\E \h
\mbox{ in $L^1(Q_{T})$}.\]

Finally, we show that $(\E,\B)$ is a weak solution of (1.1)-(1.4).
By multiplying Eq.(2.11) and Eq.(2.12) by 
test functions $\Ps$ and $\Ph$, respectively,
and then taking integration over
$Q_{T}$, after some routine calculations and taking the limit, we see 
that $(\E,\B)$ 
is a weak solution to the system (1.1)-(1.4). 

Now we consider the case where $\sigma(x,s)$ is discontinuous on $s$ 
at some points.
Without loss of generality, we may assume that $\sigma(x,s)$ has a jump only
at one point $s=1$. In this case $\sigma(x,s)$ is not uniquely defined at
$s=1$. We shall understand the value of $\sigma(x,1)$ in the following sense: 
\[ \sigma(x,1)\in [\sigma(x,1-),\sigma(x,1+)],\]
where $\sigma(x, 1\pm)$ represents the right or left limit as $s\rightarrow 1$.

By the standard approximation, we can construct a smooth 
approximation sequence $\sigma_m(x,s)$ such that 
\bys
(i) & & \sigma_m(x,s)\, \mbox{ is monotonic increasing for all $s\geq 0$},\\
(ii) & & \sigma_m(x,s)=\sigma(x,s), \mbox{ if $|s-1|\geq \frac{1}{m}$}.
\eys

Let $(\E_m,\vh_m)$ be a solution of (1.1)-(1.4) in which $\sigma(x,s)$ is 
replaced by 
$\sigma_m(x,s)$. By the same energy estimate we see that there exists a 
measurable function $\beta(x,t)\in L^{\frac{p+2}{p+1}}(Q_{T})$ such that
\[ \sigma_m(x,|\E_m|)\rightarrow \beta(x,t), \mbox{weakly in 
$L^{\frac{p+2}{p+1}}(Q_{T})$}.\]
Define
\bys
& & A_m=\{ (x,t): 1-\frac{1}{m}\leq |\E|\leq 1+\frac{1}{m}\},\\
& & A=\{ (x,t): |\E(x,t)|=1\}.
\eys
Since $\sigma(x,s)$ is continuous except at $s=1$, we see
\[ \beta(x,t)=\sigma(x, |\E|)\E, \mbox{if $(x,t)\in Q_{T}\backslash A$}.\]
Now it is clear
that \[ A=\bigcap_{m=1}^{\infty} A_m.\]
Recall that $\sigma_m(x,s)=\sigma(x,s)$ if $|s-1|\geq \frac{1}{m}$. 
It follows that for all $(x,t)\in A_m$ 
\[ \sigma(x, 1-\frac{1}{m})\leq \sigma_m(x,|\E|)\leq \sigma(x,1+\frac{1}{m}).\]
Consequently, as $m\rightarrow 0$,
\[ \sigma(x, 1-)\leq \beta(x,t)\leq \sigma(x,1+), (x,t)\in A.\]
Thus, $(\E, \vh)$ is a weak solution of (1.1)-(1.4).

Finally, we show the uniqueness. Suppose $(\E, \vh)$ and $(\E^*, \vh^*)$
are two solutions of (1.1)-(1.4). Let
\[ \hat{\E}=\E-\E^*, \hat{\vh}=\vh-\vh^*.\]
Similar to the calculation in deriving energy estimates, we find
\bys
&& \sup_{0\leq t\leq T}\int_{\Omega}
[|\hat{\E}|^2+|\hat{\vh}|^2] dx+\int\int_{Q_{T}}\left[
\sigma(x, |\E|)\E-\sigma(x, |\E^*|)\E^*\right]\cdot
\left[\E-\E^*\right] dxdt\\
& & \leq 0.
\eys
The monotonicity of $\sigma(x, s)$ implies that the second term in the above
inequality is nonnegative. It follows that
\[ \sup_{0\leq t\leq T}\int_{\Omega}
[|\hat{\E}|^2+|\hat{\vh}|^2] dx\leq 0.\]
Therefore, the uniqueness follows immediately.

\hfill Q.E.D.

\begin{center}
{\bf 3. Singular Limit Problem}
\end{center}

In this section we shall show that the solution
of (1.1)-(1.4) has a limit as $\z \rightarrow 0$, which
solves Maxwell's equations in quasi-stationary fields, i.e. the system
(1.1)-(1.4) with $\z=0$. A weak solution of the quasi-stationary system
is defined as in Definition 2.1 with $\z=0$.

From now on we denote by $(\E_{\z}, \B_{\z})$ the weak solution of the system
(1.1)-(1.4).

\ \\
{\bf Theorem 3.1}: {\em The limit of $(\E_{\z}, \B_{\z})$ as $\z\rightarrow 0$ 
solves the quasi-stationary
system (1.1)-(1.4) with $\z=0$ in the weak sense. Moreover, the weak solution
is unique if $\sigma(x,s)>0$ for all $(x,s)\in \Omega\times R^+$.}\\ 
{\bf Proof}: 
The
crucial step in proving the convergence is to show
\[ \sigma(|\E_{\z}|)\E_{\z}\rightarrow \sigma(|\E|)\E,\, a.e.\,
\mbox{in $Q_{T}$ 
as $\z\rightarrow
0$}.\]
The monotonicity of $\sigma(x,s)$ in $s$ plays a key rule.

First of all, from Lemma 2.1 and the weak compactness we see
\bys
& & \E_{\z}\rightarrow \E,\, \B_{\z t}\rightarrow \B_{t}, 
\mbox{weakly in $L^2(Q_{T})$},\\
& & \g \B_{\z} \rightarrow \g \B,\, 
\mbox{weakly in $L^{\frac{p+2}{p+1}}(Q_{T})$},\\
& & \sigma(x,|\E_{\z}|)\E_{\z}\rightarrow {\bf J}(x,t), 
\h \mbox{ weakly in $L^{
\frac{p+2}{p+1}}(Q_{T})$},
\eys
where ${\bf J}(x,t)\in L^{\frac{p+2}{p+1}}(Q_{T})$.
Moreover, as $div \B_{\z}(x,t)=0$, by the decomposition property 
of $H^1(\Omega)$ property,
after extracting a subsequence if necessary we see that
\[ \B_{\z}\rightarrow \B,\, \mbox{strongly in $L^2(Q_{T})$ }\]
and
\[ \B_{\z}\rightarrow \B, \, a.e.\, \mbox{ in $Q_{T}$}.\]

Next we show 
\[ {\bf J}(x,t)=\sigma(x,|\E|)\E, \h a.e. \mbox{in} \, Q_{T}.\]

We use a monotonicity argument. 
As a first step, we show
\[ \lim_{\z\rightarrow 0}\int\int_{Q_{T}}\sigma(x,|\E_{\z}|)|\E_{\z}|^2dxdt
=\int\int_{Q_{T}}\J\cdot \E dxdt.\]
Here we adopt an idea from \cite{JO}.
Let $\lambda(t)$ be a nonnegative smooth function and 
\[ \lambda'(t)\leq 0,
\lambda(0)=1, \lambda(T)=0.\]
Define an operator ${\cal}{L}$ in $L^{p+2}(Q_{T})^3$ as follows:
\[ {\cal}{L} [\E]=\sigma(x,|\E|)\E .\]
Since $\sigma(x, s)$ is monotonic increasing in $s$, 
then the operator ${\cal}{L}$ is monotonic increasing, that is, 
\bys
& & <{\cal}{L}[\E_{\z}]-{\cal}{L}[\E], \E_{\z}-\E>\geq 0.
\eys

It is clear that
\bys
& & <{\cal}{L}[\E_{\z}]-{\cal}{L}[\E], \E_{\z}-\E>\\
& & =< {\cal}{L}[\E_{\z}], \E_{\z}>-<{\cal}{L}[\E_{\z}], \E>-<{\cal}{L}
[\E],\E_{\z}>+<{\cal}{L}[\E],\E>\\
\eys
It follows that 
\setcounter{section}{3}
\setcounter{local}{1}
\by
& & \lim_{\z\rightarrow 0}inf <{\cal}{L}[\E_{\z}], \E_{\z}>
\geq <{\bf J},
\E>.
\ey
On the other hand,
from the system (1.1)-(1.2) we have
\by
& & \int_{0}^{T}\int_{\Omega} \sigma(\E_{\z})|\E_{\z}|^2\lambda(t)dxdt\nonumber
\\
& & =- \int_{0}^{T}\int_{\Omega}\left[
\z \lambda(t)\E_{\z t}\cdot \E_{\z}+\lambda \vh_{\z t}\cdot \vh_{\z}\right]dxdt
+\int\int_{Q_{T}}\F \cdot \E_{\z}dxdt\nonumber\\
& & =-\int_{0}^{T}\int_{\Omega}
\left\{ \frac{\partial}{\partial t}\left[ \frac{1}{2}(\z |\E_{\z}|^2+|\vh_{\z}|^2)\lambda(t)\right]-\lambda'(t)\left[ \frac{1}{2}(\z |\E_{\z}|^2+|\vh_{\z}|^2)
\right]\right\}dxdt+\nonumber\\
& & \int\int_{Q_{T}}\F \cdot \E_{\z}dxdt \nonumber \\
& & \leq \int_{0}^{T}\int_{\Omega}
\left\{ \lambda'(t)\left[ \frac{1}{2}|\vh_{\z}|^2
\right]\right\}dxdt+\frac{1}{2} \int_{\Omega}\left[ \z |\E_{\z}(x,0)|^2+|\vh_{\z}(x,0)|^2\right] dx+\nonumber\\
& & \int\int_{Q_{T}}\F \cdot \E_{\z}dxdt.\scl  
\ey
Since $\lambda'(t)\leq 0$, it follows that
\by
& &  \lim_{\z\rightarrow 0}sup \int_{0}^{T}\int_{\Omega} \sigma(|\E_{\z}|)|\E_{\z}|^2
\lambda'(t)dxdt\nonumber\\
& & \leq \frac{1}{2}\int_{0}^{T}\int_{\Omega}
\lambda(t)|\vh|^2
dxdt+\frac{1}{2} \int_{\Omega}|\vh_{0}(x)|^2dx+\int\int_{Q_{T}}
\E\cdot \F dxdt.\scl
\ey

Recall from Definition 2.1 
that $(\E_\z, \B_{\z})$ satisfies the following integral equations:
\bys
& & \int_{0}^{T}
\int_{\Omega}\left[-\z \E_{\z}\cdot {\bf \Phi}_{t}+ \sigma(x,|\E_{\z}|)\E_{\z}
\cdot
\Ph \right] dxdt\\
& & = \int_{0}^{T}\int_{\Omega}{\bf H_{\z}}\cdot (\g \Ph)dxdt
+\int\int_{Q_{T}}\F \cdot \Ph_{\z}dxdt+\int_{\Omega}\z \E_0\cdot \Ph(x,0)dx,\\
& & \int_{0}^{T}\int_{\Omega}\left[- \B_{\z}
\cdot \Ps_{t} + \E_{\z} \cdot (\g \Ps)\right] dxdt=
\int_{\Omega} \left[ {\bf H}_{0}(x) \cdot \Ps(x,0)\right]dx. 
\eys
Note by Lemma 2.1 that
\bys
& & |\int_{0}^{T}\int_{\Omega} \z \E_{\z t}\cdot \Ph dxdt|\\
& & \leq \left(\int_{0}^{T}\int_{\Omega} \z |\E_{\z t}|^2dxdt\right)^{1/2}
\left(\int_{0}^{T}\int_{\Omega} \z |\Ph|^2 dxdt\right)^{1/2}\\
& & \rightarrow 0 \h \mbox{as $\z\rightarrow 0$}.
\eys

We take the limit as $\z\rightarrow 0$ to obtain
\by
& & \int_{0}^{T}
\int_{\Omega}\left[{\bf J}
\cdot
\Ph \right] dxdt = \int_{0}^{T}\int_{\Omega}{\bf H}\cdot (\g \Ph)dxdt+
\int\int_{Q_{T}}\F \cdot \Ph dxdt\scl\\
& & \int_{0}^{T}\int_{\Omega}\left[- \B
\cdot \Ps_{t} + \E\cdot (\g \Ps)\right] dxdt=
\int_{\Omega} \left[ {\bf H}_{0}(x) \cdot \Ps(x,0)\right]dx. \scl
\ey
Now by choosing $\Ph =\lambda(t) \E$ and $\Ps=\lambda(t) \vh$
(note that the condition at $t=T$ is satisfied since $\lambda(T)=0$),
we obtain
\[ 
\int_{0}^{T}\int_{\Omega}[\lambda(t) {\bf J}\cdot \E] dxdt
=\frac{1}{2}\int_{0}^{T}\int_{\Omega}
[ \lambda'(t)|\vh|^2]
dxdt+\frac{1}{2} \int_{\Omega}|\vh_{0}(x)|^2] dx+\int\int_{Q_{T}}
\E\cdot\F dxdt.\]
It follows by (3.3) that 
\by
\lim_{\z \rightarrow 0}sup <\lambda(t){\cal}{L}[\E_{\z}], \E_{\z}> 
\leq <\lambda(t){\bf J}, \E>.\scl
\ey
Combination of (3.1) and (3.6) yields
\by
\lim_{\z \rightarrow 0}<\lambda(t){\cal}{L}[\E_{\z}], \E_{\z}>
=<\lambda(t){\bf J}, \E>.\scl
\ey

Consequently, by choosing $\lambda(t)$ properly we have
\by
& & \lim_{\z \rightarrow 0}<{\cal}{L}[\E_{\z}], \E_{\z}>
=<{\bf J}, \E>.\scl
\ey


For any vector field $\W\in L^2(Q_{T})\bigcap L^{p+1}(Q_{T})$,
the monotonicity of $\sigma(x,s)$ in $s$ implies 
\by
\int\int_{Q_{T}}[\sigma(x, |\E_{\z}|)\E_{\z}-\sigma(x, |\W|)\W]
\cdot [\E_n-\W]dxdt\geq 0.\scl
\ey
We rewrite the above inequality to the following form:
\by
& & \int\int_{Q_{T}}\left\{
\sigma(|\E_{\z}|)|\E_{\z}|^2-\sigma(|\E_{\z}|)\E_{\z} \cdot \W\right\}dxdt
\nonumber\\
& & \geq \int\int_{Q_{T}}[\sigma(x,|\W|)\W]\cdot [
\E_{\z}-\W]dxdt.\scl
\ey
We take the limit as $\z\rightarrow 0$ and use (3.8) for the first term in 
(3.10) to obtain
\bys
& & \int\int_{Q_{T}}\left\{
{\bf J}\cdot \E-{\bf J}\cdot \W\right\}dxdt\nonumber\\
& & \geq \int\int_{Q_{T}}[\sigma(x, |\W|)\W]\cdot [
\E-\W]dxdt.
\eys
Equivalently,
\by
& & \int\int_{Q_{T}}\left\{
[{\bf J}-
\sigma(x, |\W|)\W]\cdot [
\E-\W]\right\} dxdt\geq 0.\scl
\ey
Set $\W=\E+\delta \Y $, 
where $\delta>0$ is small parameter and $\Y \in
L^2(Q_{T})\bigcap L^{p+1}(Q_{T})$ is arbitrary.

With the above choice of $\W$ in the equality (3.10), we obtain
\bys
& & \int\int_{Q_{T}}\left\{
\Y \cdot [{\bf J}-\sigma(|\E+\delta \Y |)[\E+\delta \Y ]
\right\} dxdt\geq 0.
\eys
When $\sigma(x,s)$ is continuous in $s$, then
we let $\delta\rightarrow 0$ to obtain
\[ {\bf J}(x,t)=\sigma(x,|\E|)\E,\]
since $\Y(x,t)$ is arbitrary in $L^2(Q_{T})\bigcap L^{p+1}(Q_{T})$.

When $\sigma(x,s)$ has a jump at a point, say, $s=1$. Then as in \S 2 we understand
that the value of $\sigma(x, s)$ at $s=1$ is
\[ \sigma(x,1-)\leq \sigma(x,1)\leq \sigma(x,1+).\]
By using the same procedure as in \S 2, we can derive the above inequality.

Finally, by taking limit for (2.1)-(2.2) we see that
$(\E, \vh)$ is a weak solution of the quasi-stationary system. 

To prove the uniqueness, we assume 
that $(\E,\B)$ and $(\E^*,\B^*)$ are two weak solutions to
the quasi-stationary system. Let $\hat{\E}=\E-\E^*$ and $\hat{\B}=
\B-\B^*$. Then the energy estimate implies
\bys
& & \sup_{0\leq t\leq T}\int_{\Omega}|\hat{\B}|^2dx+\int\int_{Q_{T}}
[\sigma(x,|\E|)\E-\sigma(x, |\E^*|)\E^*]\cdot [\E-\E^*]dxdt\leq 0.
\eys
The monotonicity of $\sigma(x,s)$ in $s$ implies the second term in
the above inequality is nonnegative. It follows that
\[ \hat{\B}=0,\h a.e. \mbox{in $Q_{T}$}.\]
From the definition of weak solution, we have 
$\hat{\E}=0$ as long as $\sigma>0$.

\hfill Q.E.D.

 To conclude this section, 
 we consider two special classes of electric conductivity
$ \sigma(x, s)=s^{p}$ and 
\[ 
\sigma(x,s)=\left\{ \begin{array}{ll}
& a, \,\mbox{if $|s|\leq 1$};\\
& b, \,\mbox{if $|s|>1$},
\end{array}
\right.
\]
where $0<a<b$.

For the first case with $\F=0$, it is easy to see that
\[ \E=|\g \B|^{-\frac{p}{p+1}}\g\B.\]
It follows from Theorem 3.1 that there exists a unique weak solution to
the following evolution system
\[ \B_{t}+\g[ |\g \B|^{-\frac{p}{p+1}}\g\B]={\bf G}, (x,t)\in Q_{T},\]
subject to the initial-boundary conditions:
\bys
& & {\bf N}\times (\g\B)=0,\hh (x,t)\in \partial \Omega\times [0,T],\\
& & \B(x,0)=\B_0(x),\hh x\in \Omega,
\eys
where ${\bf G}$ is a known exterior magnetic field. This existence result is
not covered in \cite{Y2}. More regularity of the weak solution can be
established as in \cite{Y2}. We shall not repeat it here.

For the second case, we see that
\[ \E= \frac{1}{a} \g\B+\frac{1}{a}\F,\]
in the region $ Q_{T}^{-}=\{ (x,t): |\E|<1\}$ and
\[ \E =\frac{1}{b} \g\B+\frac{1}{b}\F,\]
in the region $Q_{T}^{+}=\{ (x,t): |\E|>1\}$.
Note that
\[ div\B(x,t)=div\B_0(x)=0,\]
we see that
\[ \g\g\B=-\Delta \B.\]
It follows that $\B$ satisfies the parabolic equation:
\[ \B_{t}-\frac{1}{a}\Delta \B=-\frac{1}{a}\g\F, \h (x,t)\in Q_{T}^{-}\]
and
\[  \B_{t}-\frac{1}{b}\Delta \B=-\frac{1}{b}\g\F, \h (x,t)\in Q_{T}^{+}.\]
The regularity theory of parabolic equations implies that
$\B(x,t)$ is smooth in $Q_{T}^{\pm}$.

The interface between $Q_{T}^{-}$ and $Q_{T}^{+}$ is defined by
\[ \Gamma=\{(x,t): |\E|=1\},\]
which is a free boundary. 

\ \\
{\bf Remark 3.1}: {\em We may allow that $\sigma(s)=0$ if $|\E|<1$ in the
above example. In this case, one must consider the full system (1.1)-(1.2)
in order to define a weak solution. However, in this case 
the uniqueness of the
weak solution does not hold.}
\ \\
{\bf Remark 3.2}: 
{\em It is not clear whether or not $\Gamma$ is indeed a hypersurface
in $R^3\times (0,\infty]$.
It would be of great interesting to study the
smoothness of the interface $\Gamma$ and to find the free boundary conditions
for $\B$.}

\ \\
{\bf Acknowledgment}: {\em This work was done when author was visiting
University of California and MSRI at Berkeley.
The author would like to thank them for the hospitality. The author
is also grateful to Professor L.C. Evans and Professor Guy David for
many helpful discussions.}

\end{document}